\documentclass[oneside,english]{article}
\usepackage[T1]{fontenc}
\usepackage{amssymb,amsthm,amsmath}
\usepackage{color}
\usepackage{babel}
\usepackage[dvips]{hyperref}

\title{Occupation time fluctuations of Poisson and equilibrium branching systems in critical and large dimensions}
\author{Piotr Mi\l{}o\'s \\
Institute of Mathematics\\
Polish Academy of Sciences\\Warsaw}
 \theoremstyle{plain}
 \newtheorem{thm}{Theorem}[section]
 \numberwithin{equation}{section}
 \numberwithin{figure}{section}
 \theoremstyle{plain}
 \theoremstyle{remark}
 \newtheorem{rem}[thm]{Remark}
 \theoremstyle{plain}
 \newtheorem{fact}[thm]{Fact}
 \theoremstyle{plain}
 
 \theoremstyle{remark}
 \newtheorem*{acknowledgement*}{Acknowledgement}

\newcommand{\ev}[1]{\mathbb{E}{#1}}

\newcommand{\norm}[3]{\Vert #1 \Vert_{ #2 } ^{ #3 }}
\newcommand{\rbr}[1]{\left( #1 \right)}
\newcommand{\ddp}[2]{\left\langle #1, #2 \right\rangle}
\newcommand{\intr}{\int_{\mathbb{R}^d}}
\newcommand{\inti}{\int_{0}^{+\infty}}
\newcommand{\intc}[1]{\int_{0}^{#1}}
\newcommand{\Rd}{\mathbb{R}^d}
\newcommand{\T}[1]{\mathcal{T}_{#1}}

\newcommand{\norma}[1]{\left| #1 \right|^\alpha}

\newcommand{\TODO}[1]{\bgroup
    \definecolor{colabc}{rgb}{1,0,0}
    \color{colabc}
    TODO - #1
\egroup
}

\newcommand{\Tphi}{\varphi_T}
\newcommand{\Tchi}{\chi_T}

\begin{document}
\maketitle
\begin{abstract}
Limit theorems are presented for the rescaled occupation time fluctuation process of a critical finite variance branching particle system in $\mathbb{R}^{d}$ with symmetric $\alpha$-stable motion starting off from either a standard Poisson random field or the equilibrium distribution for critical $d=2\alpha$  and large $d>2\alpha$ dimensions. The limit processes are generalised Wiener processes. The obtained convergence is in  space-time, finite-dimensional distributions sense. With the addtional assumption on the branching law we obtain functional convergence.
\end{abstract}



\section{Introduction}
The basic object of our investigation is a branching particle system. It consists of particles evolving independently in $\Rd$ according to a spherically symmetric $\alpha$-stable Lévy process (called standard $\alpha$-stable process), $0<\alpha\leq 2$. The system starts off at time $0$ from a random point measure $M$. The lifetime of a particle is an exponential random variable with parameter $V$. After that time the particle splits according to the law determined by a generating function $F$. We always assume that the branching is critical, i.e., $F'(0)=1$. Each of the new-born particles undertake the $\alpha$-stable movement independently of the others, and so on. The evolution of the system is described by (and in fact can be identified with) the empirical (measure-valued) process $N$, where $N_t(A)$ denotes the number of particles in the set $A\subset\Rd$ at time $t$. We define the rescaled occupation time fluctuation process by
\begin{equation}
X_T(t) = \frac{1}{F_T} \intc{Tt} \left( N_s - \ev N_s \right) ds, \: t \geq 0, \label{def:occupation}
\end{equation}
where $T$ is a scaling parameter which accelerates time ($T\rightarrow +\infty$) and $F_T$ is a proper deterministic norming. $X_T$ is a signed-measure-valued process but it is convenient to regard it as a process in the tempered distributions space $S'(\Rd)$. The objectives are to find suitable $F_T$ such that $X_T$ converges in law as $T\rightarrow +\infty$ to a non-trivial limit and to identify this limit. This problem, or its modifications (e.g. its superprocess or discrete versions) has been studied in several papers (\cite{DW}, \cite{BZ2}, \cite{M}, the list is not complete). Typically, the initial configuration $M$ was a Poisson measure, in most cases a homogeneous one, i.e. with the intensity measure $\lambda$=Lebesgue measure, and the branching law was either binary or of a special form, belonging to the domain of attraction of a $(1+\beta)$ stable distribution ($0<\beta\leq 1$).
We consider a general branching law with finite variance and the initial measure $M$ is either Poisson homogeneous or is the equilibrium measure of the system. In what follows, we will use superscripts $Poiss$ (e.g., $N^{Poiss}$) or $eq$ (e.g., $X_T^{eq}$) to indicate which model we are dealing with.\\
It is known \cite{GW2} that an equilibrium measure $M^{eq}$ of our branching system exists provided that $d>\alpha$. In \cite{M} the case of intermediate dimensions $\alpha < d < 2 \alpha$ was considered. It was shown that the limits (in the sense of the convergence in law in $C([0,\tau], S'(\Rd)), \tau>0$ of $X_T^{Poiss}$ and $X_T^{eq}$ are different; they have the form $K\lambda\xi$, where $K$ is a constant and $\xi$ is a real Gaussian process which in the Poisson case is a sub-fractional Brownian motion, while in the equilibrium case it is a fractional Brownian motion (see \cite{BGT1} for the definition and properties of the sub-fractional motion). \\
In the present paper we consider the case of critical ($d=2\alpha$) and large ($d>2\alpha$) dimensions. It turns out that now the limits of $X^{Poiss}$ and $X^{eq}$ coincide, for $d=2\alpha$ the limit is $K\lambda\beta$ where $\beta$ is the standard Brownian motion and if $d>2\alpha$ then the limit is an $S'(\Rd)$-valued Wiener process. Moreover, these limits are, up to a constant, the same as those obtained in \cite{BGT3} for the Poisson system with binary branching. The method of the proof is similar to that employed in \cite{BGT3}, based on the so-called space-time approach, but in the present case some new technical difficulties appear, especially in the study of the system in equilibrium.
\section{Results}
As mentioned in Introduction our basic state space is the space $S'(\Rd)$ of tempered distributions, dual to the space $S(\Rd)$ of smooth rapidly decreasing functions. Duality in the appropriate spaces is denoted by $\ddp{\cdot}{\cdot}$. Three kinds of convergence are used. Firstly the convergence of finite-dimensional distributions, denoted by $\Rightarrow_f$.
For a continuous, $S'(\Rd)$-valued process $X = (X_t)_{t\geq 0}$ and any $\tau>0$ one can define an $\mathcal{S}'(\mathbb{R}^{d+1})$-valued random variable
\begin{equation}
\ddp{\tilde{X}}{\Phi} = \intc{\tau} \ddp{X_s}{\Phi(\cdot, s)} ds, \: \Phi \in \mathcal{S}(\mathbb{R}^{d+1}) \label{def: space-time}
\end{equation}
If for any $\tau>0$ $\tilde{X}_n\rightarrow \tilde{X}$ in distribution then we say that the convergence  in the space-time sense holds and denote this fact by $\Rightarrow_i$. Finally, we consider the functional weak convergence denoted by $X_n\Rightarrow_c X$. It holds if for any $\tau>0$ processes $X_n = (X_n(t))_{t\in [0, \tau]}$ converge to $X = (X(t))_{t\in [0, \tau]}$ weakly in $C([0, \tau], \mathcal{S}'(\Rd))$. 
It is known that $\Rightarrow_i$ and $\Rightarrow_f$ do not imply each other, but either of them together with tightness implies $\Rightarrow_c$ (\cite{BGR}). Conversely, $\Rightarrow_c$ implies both $\Rightarrow_i$, $\Rightarrow_f$. Consider a branching particle system described in Introduction.\\ Denote (recall that $F$ is the generating function of the branching law)
\begin{equation}
 m=F''(1) \label{def: m}
\end{equation}
We start with the large dimension case.
\begin{thm} \label{thm: big-dimensions}
 Assume $d>2\alpha$ and let $F_T = T^{\frac{1}{2}}$. Assume that the initial configuration of the system is given either by a Poisson homogeneous measure or by the equilibrium measure and let $X_T$ be defined by (\ref{def:occupation}), i.e. $X_T = X_T^{Poiss}$ or $X_T = X_T^{eq}$. Then
\begin{enumerate}
\item $X_T \Rightarrow_f X \text{ and } X_T \Rightarrow_i X \text{ as }  T\rightarrow +\infty$\\
where $X$ is a centered $\mathcal{S}'$-valued, Gaussian process with the covariance function
\begin{equation}
	Cov\left( \ddp{X_s}{\varphi_1}, \ddp{X_t}{\varphi_2} \right) = 
(s\wedge t) \frac{1}{2\pi} 
\intr \left( \frac{2}{|z|^\alpha} + \frac{Vm}{2|z|^{2\alpha}} \right) \widehat{\varphi_1}(z) \overline{\widehat{\varphi_2}(z)} dz \nonumber
\end{equation}
where$\ \varphi_1,\varphi_2 \in\mathcal{S}\left(\mathbb{R}^{d}\right)$.
\item If, additionally, the branching law has finite the fourth moment then
\[X_T \Rightarrow_c X \text{ as }  T\rightarrow +\infty\]
\end{enumerate}
\end{thm}
For the critical dimension we have the following theorem
\begin{thm} \label{thm: critical-dimensions}
 Assume $d=2\alpha$ and let $F_T = (T\log T)^{\frac{1}{2}}$. Assume that the initial configuration of the system is given either by a Poisson homogeneous measure or by the equilibrium measure and let $X_T$ be defined by (\ref{def:occupation}), i.e. $X_T = X_T^{Poiss}$ or $X_T = X_T^{eq}$. Then
\begin{enumerate}
\item $X_T \Rightarrow_f X \text{ and } X_T \Rightarrow_i X \text{ as }  T\rightarrow +\infty$\\
where 
\[X =  \rbr{\frac{mV}{2}}^{1/2} C_d \lambda \beta \] 
\[C_d = \left( 2^{d-2} \pi^{d/2} d \Gamma \left( \frac{d}{2} \right)\right)^{-1/2}\]
and $\beta$ is a standard Brownian motion.
\item If, additionally, the branching law has finite the fourth moment then
\[X_T \Rightarrow_c X \text{ as }  T\rightarrow +\infty\]
\end{enumerate}
\end{thm}
\begin{rem}
\begin{itemize}
\item[(a)] It is unclear if the assumption of the existence of the fourth moment is necessary for the functional convergence in both theorems to hold. One can see that only the second moment influences the result. In the proof below the assumption is only used in the proof of tightness of family $X_T$. (see also Remark \ref{rem:fourth}).
\item[(b)] The limit process $X$ in Theorem 2.1 is an $S'(\Rd)$-valued homogeneous Wiener process.
\end{itemize}
\end{rem}

\section{Proofs}
\subsection{General Scheme}
\subsubsection{Space-time convergence}\label{sec:spacetime-schema}
We present a general scheme which will be used in the proofs of both theorems. It is similar to the one employed in \cite{M} and \cite{BGT3}. Many parts of the proofs are the same for $N^{Poiss}$-system starting from a Poisson field and $N^{eq}$-system starting from the equilibrium distribution hence we will omit superscripts when a formula holds for both of them. Let $X_T$ be the occupation time fluctuation process defined by (\ref{def:occupation}). Firstly we establish the convergence in the space-time sense. Let us consider $\tilde{X}_T$ defined according to (\ref{def: space-time}). We show the convergence of the Laplace transforms
\begin{equation}
\lim_{T\rightarrow + \infty} \ev{e^{-\ddp{\tilde{X}_T}{\Phi}}} = \ev{e^{-\ddp{\tilde{X}}{\Phi}}},\: \Phi\in \mathcal{S}(\mathbb{R}^{d+1}), \: \Phi \geq 0, \label{def: laplace-convergence}
\end{equation}
where $X$ is the corresponding limit process. This will imply the weak convergence of $\tilde{X}_T$ since the limit processes are Gaussian ones (see detailed explanation in \cite{BGT2}). The purpose of the rest of this section is to gather facts used to calculate the Laplace transforms and to show convergence (\ref{def: laplace-convergence}).
To make the proof shorter we will consider $\Phi$ of the special form:
\[
\Phi(x,t) = \varphi(x) \psi(t)\:\:\varphi\in \mathcal{S}(\Rd), \psi\in \mathcal{S}(\mathbb{R}^+), \varphi\geq 0, \phi \geq 0.
\]
We also denote
\begin{equation}
\varphi_T = \frac{1}{F_T}\varphi, \: \chi(t) = \int_t^1 \psi(s) ds, \: \chi_T(t) = \chi(\frac{t}{T}). \label{def: abbrev} 
\end{equation}
We write
\begin{equation}
\Psi(x,t) = \varphi(x) \chi(t), \label{def: Psi_bez_T}
\end{equation}
\begin{equation}
\Psi_T(x,t) = \varphi_T(x) \chi_T(t), \label{def: Psi}
\end{equation}
note that $\Psi$ and $\Psi_T$ are positive functions.
For a generating function $F$ we define
\begin{equation}
G\left(s\right)=F\left(1-s\right)-1+s.\label{def: G}
\end{equation}
We will need following properties of $G$ (we omit proofs which are straightforward)
\begin{fact}
\label{fact: properties of G}
\begin{enumerate}
\item $G\left(0\right)=F\left(1\right)-1=0$,
\item $G'\left(0\right)=-F'\left(1\right)+1=0$
\item $G''\left(0\right)=F''\left(1\right)<+\infty$,
\item $G\left(v\right)=\frac{m}{2}v^{2}+g\left(v\right)v^{2}$ where $m$ is defined by (\ref{def: m}) and $\lim_{v\rightarrow0}g\left(v\right)=0$.
\item $G'''(0)<+\infty$ and $G^{IV}(0)<+\infty$ if the law determined by $F$ has finite the fourth moment
\end{enumerate}
\end{fact}
Let us recall classical Young's inequality
\begin{equation}
\norm{f \ast g}{p}{} \leq \norm{f}{q_1}{} \norm{g}{q_2}{},\:\:  \label{ineq:young}
\end{equation}
which holds when $\frac{1}{p} = \frac{1}{q_1} + \frac{1}{q_2} - 1, \: q_1, q_2 \geq 1$.\\
Now we introduce an important function used throughout the rest of the paper
\begin{equation}
 v_{\Psi}\left(x,r,t\right)=1-\mathbb{E}\exp\left\{ -\int_{0}^{t}\left\langle N_{s}^{x},\Psi\left(\cdot,r+s\right)\right\rangle ds\right\} \label{def: v}
\end{equation}
where $N_{s}^{x}$ denotes the empirical measure of the particle system with the initial condition $N_{0}^{x}=\delta_{x}$. $v_{\Psi}$ satisfies the equation
\begin{equation}
v_{\Psi}\left(x,r,t\right)=\int_{0}^{t}\mathcal{T}_{t-s}\left[\Psi\left(\cdot,r+t-s\right)\left(1-v_{\Psi}\left(\cdot,r+t-s,s\right)\right)-VG\left(v_{\Psi}\left(\cdot,r+t-s,s\right)\right)\right]\left(x\right)ds.\label{eq: v}
\end{equation}
The equation can be proved using the Feynman-Kac formula in the same way as \cite[Lemma 3.4]{M}. We also define
\begin{equation}
n_{\Psi}\left(x,r,t\right)=\int_{0}^{t}\mathcal{T}_{t-s}\Psi\left(\cdot,r+t-s\right)\left(x\right)ds.\label{def: n}\end{equation}
Since we consider only positive $\Psi$ hence (\ref{def: v}) and (\ref{eq: v}) yield
\begin{equation}
0\leq v_T(x,r,t) \leq n_T(x,r,t) \label{ineq: n<v} 
\end{equation}
where, for simplicity of notation, we write \begin{equation}
v_{T}\left(x,r,t\right)=v_{\Psi_{T}}\left(x,r,t\right),\label{eq:v_T notacja}\end{equation}
 \begin{equation}
n_{T}\left(x,r,t\right)=n_{\Psi_{T}}\left(x,r,t\right),\label{eq:n_T notacja}\end{equation}
 \begin{equation}
v_{T}\left(x\right)=v_{T}\left(x,0,T\right),\label{eq: v_T_x notacja}\end{equation}
 \begin{equation}
n_{T}\left(x\right)=n_{T}\left(x,0,T\right)\label{eq: n_T_x notacja}\end{equation}
 when no confusion can arise.\\
\begin{fact}
$n_{T}\left(x,T-s,s\right)\rightarrow0$ uniformly in $x\in\mathbb{R}^{d}$,
$s\in\left[0,T\right]$ as $T\rightarrow+\infty$. \label{fact: uniformaly}
\end{fact}
The proof is the same as \cite[Fact 3.7]{M}.\\
We also introduce function $V_T$. It is defined by
\begin{equation}
V_T(x,l) = 1 - \ev{exp\rbr{\ddp{N_l^x}{\ln (1-v_T)}}} \label{def:V}
\end{equation}
and fulfills the equation
\begin{equation}
V_{T}\left(x,l\right)=\mathcal{T}_{l}v_{T}\left(x\right)-V\int_{0}^{l}\mathcal{T}_{l-s}G\left(V_{T}\left(\cdot,s\right)\right)\left(x\right)ds.\label{main_equation}
\end{equation}
It satisfies (details can be found in \cite[Section 3.2.2]{M})
\begin{equation}
	0\leq V_{T}\left(x,l\right)\leq\mathcal{T}_{l}v_{T}\left(x\right),\,\forall_{x\in\mathbb{R}^{d},l\geq0}.
	\label{ineq: VT<TvT}
\end{equation}
Now we can write the Laplace transforms (see \cite[Section 3.1.2, 3.2.2]{M} for calculations)
\begin{equation}
\ev{e^{-\ddp{\tilde{X}_T^{Poiss}}{\Phi}}} = e^{A(T)} \label{eq: laplace-Poiss}
\end{equation}
and
\begin{equation}
\ev{e^{-\ddp{\tilde{X}_T^{eq}}{\Phi}}} = e^{A(T) + B(T)} \label{eq: laplace-eq}
\end{equation}
where
\begin{equation}
A\left(T\right)= \int_{\mathbb{R}^{d}}\int_{0}^{T}\Psi_{T}\left(x,T-s\right)v_{T}\left(x,T-s,s\right)+VG\left(v_{T}\left(x,T-s,s\right)\right)dsdx, \label{def: A}
\end{equation}
\begin{equation}
B\left(T\right)= V\int_{0}^{+\infty}\int_{\mathbb{R}^{d}}G\left(V_{T}\left(x,t\right)\right)dxdt. \label{def: B}
\end{equation}
We consider the following decomposition of $A(T)$
\begin{equation}
A(T) =\exp\left\{ V\left(I_{1}\left(T\right)+I_{2}\left(T\right)\right)+I_{3}\left(T\right)\right\} ,\label{eq:Laplace_functional_decomposition}\end{equation}
where \begin{align}
I_{1}\left(T\right)= & \int_{0}^{T}\int_{\mathbb{R}^{d}}\frac{m}{2}\left(\int_{0}^{s}\mathcal{T}_{u}\Psi_{T}\left(\cdot,T+u-s\right)\left(x\right)du\right)^{2}dxds, \\
I_{2}\left(T\right)= & \int_{0}^{T}\int_{\mathbb{R}^{d}}\left[G\left(v_{T}\left(x,T-s,s\right)\right)-\frac{m}{2}\left(\int_{0}^{s}\mathcal{T}_{u}\Psi_{T}
\left(\cdot,T+u-s\right)\left(x\right)du\right)^{2}\right]dxds,\label{def: I2}\\
I_{3}\left(T\right)= & \int_{0}^{T}\int_{\mathbb{R}^{d}}\Psi_{T}\left(x,T-s\right)v_{T}\left(x,T-s,s\right)dxds.\end{align}
We claim that in the case of large dimensions ($d>2\alpha$) we have
\begin{eqnarray}
 I_1(T) &\rightarrow&\frac{m}{2(2\pi)^d}\intc{1}\intc{1} (r\wedge r')\psi(r) \psi(r') dr dr' \intr \frac{|\widehat{\varphi}(z)|^2}{|z|^{2\alpha}} dz  \label{res:I1_limit}\\
 I_2(T) &\rightarrow&0 \label{res:I2_limit}\\
 I_3(T) &\rightarrow& \frac{1}{(2\pi)^d}\intc{1}\intc{1} (r\wedge r')\psi(r) \psi(r') dr dr' \intr \frac{|\widehat{\varphi}(z)|^2}{|z|^{\alpha}} dz \label{res:I3_limit}
\end{eqnarray}
Using decomposition (\ref{eq:Laplace_functional_decomposition}) we obtain the limit of $A(T)$ and consequently the one for the Laplace transform (\ref{eq: laplace-Poiss}). This establishes the space-time convergence of the Poisson-starting system $X_T^{Poiss}$ considered in 1) of Theorem \ref{thm: big-dimensions}.\\
Analogously, in the critical case ($d=2\alpha$), we obtain the corresponding convergence in 1) of Theorem \ref{thm: critical-dimensions} once we show 
\begin{equation}
I_1(T) \rightarrow \frac{m}{2} C_d^2 \intc{1} \intc{1} (r\wedge r') \psi(r) \psi(r') \rbr{\intr \varphi(x) dx}^2 \label{res:I1_critical}
\end{equation}
\begin{equation}
I_2(T),I_3(T) \rightarrow 0 \label{res:I2I3_critical}
\end{equation}
The limits (\ref{res:I1_limit}) - (\ref{res:I2I3_critical}) will be obtained in Section \ref{sec:calc1} and Section \ref{sec:calc3}.\\
Now we proceed to the case of the equilibrium-starting system $X^{eq}_T$. In both Theorem \ref{thm: big-dimensions} and Theorem \ref{thm: critical-dimensions} the limits are the same as in the $X^{Poiss}_T$ case. It follows immediately from (\ref{eq: laplace-eq}) that that fact will be proved when we show
\[
B(T) \rightarrow 0\label{res:tmp_B_T}
\]
Let us first observe an elementary fact that uniform convergence $V_T(\cdot, \cdot)\rightarrow 0$ as $T\rightarrow +\infty$ holds. It is a direct consequence of Fact \ref{fact: uniformaly} and the combination of inequalities (\ref{ineq: VT<TvT}) and (\ref{ineq: n<v}). This together with Fact \ref{fact: properties of G} yields
\begin{equation}
 B(T) \leq c \inti \intr \rbr{\T{t} n_T(x)}^2 dx dt, \label{res: B_Tleq}
\end{equation}
Let us denote the right-hand side of (\ref{res: B_Tleq}) by $B_1(T)$. Now we need to obtain
\begin{equation}
 \lim_{T\rightarrow +\infty} B_1(T) = 0. \label{res:B_T_limit}
\end{equation}
which is again put off to Section \ref{sec:calc1} and Section \ref{sec:calc3}.\\
\subsubsection{Finite dimensional convergence}
A similar method, based on the Laplace transform, can be applied to prove the finite distributions convergence. Indeed, for a sequence $0\leq t_1 \leq t_2 \leq \cdots \leq t_n \leq \tau$ and functions $\varphi_1, \varphi_2,\cdots, \varphi_n \in \mathcal{S}(\Rd)$, $\varphi_i\geq0$ we write the Laplace transform
\begin{equation}
\ev{\exp\rbr{\sum_{i=1}^n \ddp{X_T(t_i)}{\varphi_i}}}. \label{def: Laplace-finite}
\end{equation}
The key observation is that, formally, 
\[
\sum_{i=1}^n \ddp{X_T(t_i)}{\varphi_i} = \ddp{\tilde{X}_T}{\Phi}
\]
if $\Phi = \sum_{i=1}^n \varphi_i \delta_{t_i}$ (which coresponds to $\Psi(x,s) = \sum_{i=1}^n \varphi_i(x) \mathbf{1}_{[0,t_i]}(s)$, recall definition (\ref{def: Psi_bez_T})). \\
It appears that the Laplace transforms (\ref{eq: laplace-Poiss}), (\ref{eq: laplace-eq}) and formulae (\ref{eq: v}), (\ref{main_equation}) are still valid for $\Phi$ and $\Psi$. The proof for the Poisson-starting system is simpler version of the one presented below and is left to the reader. We employ an approximation argument. Consider  $\Phi_n \rightarrow \Phi$ where $\Phi_n \in \mathcal{S}(\mathbb{R}^{d+1})$, additionally assume that the sequence $(\Phi_n)_n$ is chosen such that $\Psi^n(x,t) = \int_t^1 \Phi^n(x,s) ds$ is nondecreasing $\Psi^n \leq \Psi^{n+1}$. To keep the proof short we adhere to the following notation, symbols with (without) the superscript $n$ will denote functions defined for $\Phi^n$ and $\Psi^n$ (resp. $\Phi$ and $\Psi$) (eg. 
$v^n:=v_{\Psi^n}$ given by (\ref{eq: v})). $T$ is fixed and hence is omitted where possible.\\
The first assertion is that $V(x,l)$ satisfies equation (\ref{main_equation}). Definition (\ref{def:V}) implies that $V^n(x,l) \rightarrow V(x,l)$ (pointwise) which follows immediately from $v^n \rightarrow v$ (left to the reader), inequality $0\leq v \leq 1$ and the dominated convergence theorem. By assumption $\Phi^n \in \mathcal{S}(\mathbb{R}^{d+1})$ and $V^n$ satisfies the equation. Passing to the limit $n\rightarrow +\infty$ and employing the dominated convergence theorem to the right side of the equation completes the proof. \\
Now we turn to the Laplace transform (\ref{eq: laplace-eq}). It is obvious that
\[
\lim_{n} \ev{e^{-\ddp{\tilde{X}_T^{eq}}{\Phi^n}}} = \ev{e^{-\ddp{\tilde{X}_T^{eq}}{\Phi}}}.
\]
One can see that formula (\ref{eq: laplace-eq}) for $\Phi$ will be justified if only $A^n \rightarrow A$ and $B^n \rightarrow B$. To show the first one is left to the reader. It is straightforward to check that $\Phi^n\leq \Phi^{n+1}$ implies $V_n \leq V_{n+1}$ and that $G$ is nondecreasing. A standard application of the monotone convergence theorem completes the proof. The finite distributions convergence is thus established. Indeed, the above argumentation allows the calculations from Section \ref{sec:spacetime-schema} to be repeated for $\Phi = \sum_{i=1}^n \varphi_i \delta_{t_i}$ which implies the convergence of the Laplace transform (\ref{def: Laplace-finite}) and consequently the finite dimensional convergence in 1) of Theorem \ref{thm: big-dimensions} and Theorem \ref{thm: critical-dimensions}.

\subsubsection{Functional convergence} \label{sec:tightness}
In this subsection we present a general scheme of the proof of the functional convergence.
The assertion follows immediately from the part $1$ of the Theorem \ref{thm: big-dimensions} or Theorem \ref{thm: critical-dimensions}, respectively, if we prove that $\{X_T, T>2\}$ is tight in $C([0,1], \mathcal{S}'(\Rd)$ (with no loss of generality we consider $\tau = 1$). Generally, we follow the lines of the proof of tightness in Theorem 2.2 in \cite{BGT3}. However, in our case new technical difficulties arise because of a more general branching law. Some estimates are more cumbersome and some extra terms appear. Moreover, we establish tightness for $X_T^{eq}$ which was not investigated in \cite{BGT3}. This requires even more intricate computations than in the Poisson case. \\
By the Mitoma theorem \cite[Mitoma 1983]{M} it suffices to show tightness of the real processes $\left\langle X_{T},\varphi\right\rangle $
for all $\varphi\in\mathcal{S}\left(\mathbb{R}^{d}\right)$. This can be done using a criterion \cite[Theorem 12.3]{B}
\begin{equation}
 \ev{(\ddp{X_T(t)}{\varphi}, \ddp{X_T(s)}{\varphi})^4} \leq C (t-s)^2. \label{ineq:tightness}
\end{equation}
Let $(\psi_n)_n$ be a sequence in $\mathcal{S}(\mathbb{R})$, denote $\chi_n(u) = \int_u^1 \psi_n(s) ds$. It is an easy exercise to show that the sequence $(\psi_n)_n$ can be chosen in such a way that
\[
\psi_n\rightarrow \delta_t - \delta_s, 
\]
\begin{equation}
 0\leq \chi_n \leq \mathbf{1}_{[s,t]}.\label{tmp:chi_n}
\end{equation}
A detailed construction can be found in \cite{BGT3}.\\
Denote $\Phi_n=\varphi\otimes\psi_n$. We have
\[
 \lim_{n\rightarrow +\infty} \ddp{X_T}{\Phi_n} = \ddp{X_T(t)}{\varphi} - \ddp{X_T(s)}{\varphi}
\]
 thus by Fatou's lemma and the definition of $\psi_n$ we will obtain (\ref{ineq:tightness}) if we prove ($C$ is a constant independent of $n$ and $T$) that
\[
 \ev{\ddp{\tilde{X}_T}{\Phi_n}^4} \leq C (t-s)^2.
\]
From now on we fix an arbitrary $n$ and dentote $\Phi:=\Phi_n$ and $\chi := \chi_n$. 
By properties of the Laplace transform we have
\[
 \ev{\ddp{\tilde{X}_T}{\Phi}^4} = \frac{d^4}{d\theta^4}|_{\theta = 0} \ev{e^{-\theta \ddp{\tilde{X}_T}{\Phi}}}
\]
Hence the proof of tightness will be completed if we show 
\begin{equation}
 \frac{d^4}{d\theta^4}|_{\theta = 0} \ev{e^{-\theta \ddp{\tilde{X}_T}{\Phi}}} \leq C(t-s)^2, \label{eq: laplace-tighness-end}
\end{equation}
The rest of the section is devoted to calculate the fourth derivative of the Laplace transforms (\ref{eq: laplace-Poiss}) and (\ref{eq: laplace-eq}).\\
Here and subsequently $A(\theta,T)$ and $B(\theta, T)$ will denote (\ref{def: A}) and (\ref{def: B}) taken for $\Psi_{\theta,T} = \theta\varphi_T\otimes \chi_T$ ($\varphi_T$ and $\chi_T$ are defined by \ref{def: abbrev}), i.e.,
\[
A\left(\theta, T\right)=  \int_{\mathbb{R}^{d}}\int_{0}^{T}\theta \varphi_T(x) \chi_T(T-s) v_{\Psi_{\theta, T}}\left(x,T-s,s\right)+VG\left(v_{\Psi_{\theta, T}}\left(x,T-s,s\right)\right)dsdx, \label{def: Atheta}
\]
\[
B\left(\theta, T\right)= V\int_{0}^{+\infty}\int_{\mathbb{R}^{d}}G\left(V_{\Psi_{\theta, T}}\left(x,t\right)\right)dxdt. \label{def: Btheta}
\]
\begin{rem}\label{rem:fourth}
This is the point where we need the existence of the fourth moment of the branching law. Note that in the case of the binary branching law (model investigated in \cite{BGT3}) the fourth moment is obviously finite. The formulae derived below are consistent, but more complicated that the ones considered in \cite{BGT3}. This makes the computation here significantly longer and, moreover, some new technical difficulties arise especially in case of critical dimensions. New arguments and estimations were required to cope with them.
\end{rem}
A trivial verification shows that $A(0,T) = 0$, $A'(0,T) = 0$, $B(0,T) = 0$, $B'(0,T) = 0$ hence
\[
\frac{d^4}{d\theta^4} e^{A(\theta, T)} |_{\theta = 0} = A^{IV}(0,T) + A''(0,T)^2
\]
and
\[
\frac{d^4}{d\theta^4} e^{A(\theta, T) + B(\theta, T)} |_{\theta = 0} = A^{IV}(0,T) + B^{IV}(0,T) + \rbr{A''(0,T) + B''(0,T)}^2.
\]
Now taking into account (\ref{eq: laplace-tighness-end}), to demonstrate tightness, it suffices to show that
\begin{equation}
A''(0,T)\leq C(t-s), \: B''(0,T) \leq C(t-s)\label{ineq:tightness_main1}
\end{equation}
\begin{equation}
A^{IV}(0,T) \leq C (t-s)^2, \: B^{IV}(0,T) \leq C (t-s)^2\label{ineq:tightness_main2}
\end{equation}
It will convenient to denote 
\[
v(\theta) = v(\theta)(x, T-u, u) = v_{\Psi_{\theta, T}}(x, T-u,u)
\]
\[
V(\theta) = V(\theta)(x, t) =V_{\Psi_{\theta, T}}(x, T-u,u)
\]
\[
k=G'''(0), \: l=G^{IV}(0)
\]
Using the properties from Fact \ref{fact: properties of G}.
\[
A''(0,T) = 2 \intc{T} \intr \Tphi(x) \Tchi(T-u) v'(0) dx du + Vm \intc{T} \intr (v'(0))^2 dx du
\]
\begin{align*}
 A^{IV}(0,T) &= 4 \intc{T} \intr \Tphi(x) \Tchi(T-u) v'''(0) dx du + Vl \intc{T} \intr (v'(0))^4 dx du \\
	&+ 6Vk \intc{T} \intr (v'(0))^2 v''(0) dx du + 3Vm \intc{T} \intr (v''(0))^2 dx du \\
	&+ 4Vm \intc{T} \intr v'(0) v'''(0) dx du
\end{align*}
Similarly
\[
 B''(0,T) = Vm \intc{T} \intr \rbr{V'(0)}^2 ds dx
\]
\begin{eqnarray}
B^{IV}(0,T) = Vl \inti \intr \rbr{V'(0)}^4 dx ds + 6Vk \inti \intr V''(0) \rbr{V'(0)}^2 dx ds + \nonumber
\\ 3Vm \inti \intr \rbr{V''(0)}^2 dx ds + 4Vl \inti \intr V'(0) V'''(0) dx ds  \label{res:BIV}
\end{eqnarray}
Derivatives of $v(\theta)$ and $V(\theta)$ at $\theta=0$ are given by 
\begin{equation}
 v'(0)(x, T-u, u) = \intc{u} \T{u-s} [\Tphi(\cdot)\Tchi(T-s)](x) ds \label{res:v'}
\end{equation}
\begin{align*}
 v''(0)(x, T-u, u) &= -2 \intc{u} \T{u-s}[\Tphi(\cdot)\Tchi(T-s)v'(0)(\cdot,T-s,s)](x) ds \\
		&- mV \intc{u} \T{u-s}[(v'(0)(\cdot, T-s, s))^2](x) ds
\end{align*}

\begin{align*}
 v'''(0)(x, T-u, u) &= -3 \intc{u} \T{u-s}[\Tphi(\cdot)\Tchi(T-s)v''(0)(\cdot,T-s,s)](x) ds \\
		&- kV \intc{u} \T{u-s}[(v'(0)(\cdot, T-s, s))^3](x) ds \\
		&- 3mV \intc{u} \T{u-s}[v'(0)(\cdot, T-s, s) v''(0)(\cdot, T-s, s)](x) ds
\end{align*}
\[
 V'(0)(x,s) = \T{s} v'(0)(x,0,T)
\]
\begin{align}
 V''(0)(x,s) = \T{s} v''(0)(x,0,T) -Vm\intc{s} \T{t-u} \left( \rbr{V'(0)(\cdot, u)}^2 \right) du \label{res:V''}
\end{align}
\begin{align*}
 V^{IV}(0)(x,s) = \T{s} v''(0)(x,0,T) - V\intc{s} \T{t-u} \left( 3m V'(0)(\cdot, u) V''(0)(\cdot, u) +  k\rbr{V'''(0)(\cdot, u)}^3 \right) du
\end{align*}

\subsection{Proof of Theorem 2.1} \label{sec:calc1}
We follow the scheme described in Section \ref{sec:spacetime-schema} for the large dimensions case.
$I_1$ does not depend on $F$ so (\ref{res:I1_limit}) can be obtained in the same way as in \cite[(3.15)]{BGT3}.\\
We will turn now to (\ref{res:I2_limit}) which is a little more intricate.
Combining (\ref{def: I2}) and decomposition of $G$ from Fact \ref{fact: properties of G} we acquire
\begin{equation}
  I_{2}(T) =\frac{m}{2}I_{21}\left(T\right)+I_{22}\left(T\right),
\end{equation}
where 
\begin{equation}
I_{21}\left(T\right)=\int_{0}^{T}\int_{\mathbb{R}^{d}}v_{T}\left(x,T-s,s\right)^{2}-\left(\int_{0}^{s}\mathcal{T}_{u}\Psi_{T}\left(\cdot,T+u-s\right)\left(x\right)du\right)^{2}dxds, \label{eq:I21}
\end{equation}
\begin{equation}
I_{22}\left(T\right)=\int_{0}^{T}\int_{\mathbb{R}^{d}}g\left(v_{T}\left(x,T-s,s\right)\right)v_{T}\left(x,T-s,s\right)^{2}dxds\label{eq: I2pp}
\end{equation}
We have the following inequalities (proofs are straightforward and can be found \cite[Section 3.1.3]{M})
\begin{align}
 0\leq & n_{T}\left(x,T-s,s\right)-v_{T}\left(x,T-s,s\right) \leq \nonumber\\
       &C \int_{0}^{s}\mathcal{T}_{s-u}\left[\Psi_{T}\left(\cdot,T-u\right)n_{T}\left(\cdot,T-u,u\right)+n_{T}\left(\cdot,T-u,u\right)^{2}\right]\left(x\right)du \label{ineq:n-v}
\end{align}
\begin{equation}
n_{T}\left(x,T-s,s\right)+v_{T}\left(x,T-s,s\right)\leq2n_{T}\left(x,T-s,s\right)\leq2\int_{0}^{s}\mathcal{T}_{s-u}\Psi_T \left(\cdot,T-u\right)\left(x\right)du. \label{ineq:n+v}
\end{equation}
By (\ref{eq:I21}) we have
\[
 0\leq-I_{21}(T) \intc{T}\intr (n_T(x,T-s,s)-v_T(x,T-s,s))(n_T(x,T-s,s)+v_T(x,T-s,s)) ds dx
\]
We use (\ref{ineq:n-v}) and (\ref{ineq:n+v}) obtaining
\[
-I_{21}(T) \leq C(I_{211}(T) + I_{212}(T))
\]
where
\[
I_{211}(T) = \intc{T} \intr \rbr{\int_{0}^{s}\mathcal{T}_{s-u}\left[\Psi_{T}\left(\cdot,T-u\right)n_{T}\left(\cdot,T-u,u\right)\right]\left(x\right)du} \rbr{\int_{0}^{s}\mathcal{T}_{s-u}\Psi_T\left(\cdot,T-u\right)(x)du} dx ds,
\]
\[
I_{212}(T) = \intc{T} \intr \rbr{\int_{0}^{s}\mathcal{T}_{s-u}\left[n_{T}\left(\cdot,T-u,u\right)^2\right]\left(x\right)du} \rbr{\int_{0}^{s}\mathcal{T}_{s-u}\Psi_T \left(\cdot,T-u\right)(x)du} dx ds.
\]
One can see that $I_{211}$ and $I_{212}$ coincide with $J_1$ and $J_2$ from \cite{BGT3} (see (3.20) and (3.21)) hence by the proof therein
\[
\lim_{T\rightarrow + \infty} I_{21}(T) =0 
\]
Next we show that $I_{22} \rightarrow 0$. Indeed, applying Facts \ref{fact: properties of G} and \ref{fact: uniformaly} and inequality (\ref{ineq: n<v}) we obtain $\forall_{\epsilon>0} \exists_{T_0}$ such that $\forall_{T>T_0}$
\[
0\leq I_{22}(T) \leq \epsilon I_{1}(T)
\]
which clearly implies $I_{22} \rightarrow 0$.\\
Finally we obtain (\ref{res:I3_limit}). $I_3(T))$ can be split in the same way as \cite[(3.24)]{BGT3}. The only difference is that 
\[
 I_3'''(T) = \intc{T} \intr \varphi_T(x) \chi_T(T-u) \intc{u} \T{u-s} G(v_{\Psi_T}(\cdot, T-s, s))(x)ds dx du
\]
but $G(v)$ is comparable with $v^2$ hence the rest of proof goes along the same lines. (Compare to \cite[(3.27)]{BGT3}).\\
Now we turn to the equilibrium case. As observed before, it suffices to prove (\ref{res:B_T_limit}). Using the Fourier transforms we get
\[
 B_1(T) = C \int_{\mathbb{R}^{d}}\frac{1}{\left|z\right|^{\alpha}}\left(\widehat{n}_{T}\left(z\right)\right)^{2}dz
\]
It is not hard to see that
\[
 \left|\widehat{n}_{T}\left(z\right)\right|\leq \frac{C T^{1-\frac{\beta}{\alpha}}}{F_T}\frac{\left|\widehat{\varphi}\left(z\right)\right|}{|z|^\beta},\:\beta\in [0,\alpha]
\]
Hence we obtain
\[
 |B_1(T)| \leq C \frac{T^{2(1-\frac{\beta}{\alpha})}}{F_T^2} \intr \frac{|\widehat{\varphi}(z)|^2}{|z|^\alpha} \frac{1}{|z|^{2\beta}} dz.
\]
We take $\beta$ such that $\frac{1}{2}\alpha < \beta$ but $\alpha + 2\beta < d$ (it can be done because $2\alpha<d$). The first condition gives as 
\[
 \frac{T^{2(1-\frac{\beta}{\alpha})}}{F_T^2} \rightarrow 0 \qquad \text{as }T\rightarrow+\infty
\]
and the second ensures that the integral is finite. This completes the proof of (\ref{res:B_T_limit}) and consequently part 1) of Theorem \ref{thm: big-dimensions}.\\
Now we proceed to part 2). Firstly, we follow the scheme from Section \ref{sec:tightness}. The proof will be completed when we show inequalities (\ref{ineq:tightness_main1}) and (\ref{ineq:tightness_main2}). It can be done by applying the expressions derived in Section \ref{sec:tightness} repeatedly. This results in many terms which have to be estimated separately. As an example consider (\ref{res:BIV}), take only its third term then substitute $V''(0,T)$ in it utilizing only the second term of (\ref{res:V''}) and finally eliminate $v'(0,T)$ using (\ref{res:v'}). In this way we obtain
\[R=\int_{\mathbb{R}^d}\int_{0}^{+\infty} \left( \int_{0}^{l} \mathcal{T}_{l-s_1}\left[ \left( \mathcal{T}_{s_1}\left[\int_{0}^{T} \mathcal{T}_{T-s_3}\left[\varphi_{T}(\cdot)\chi_{T}(T - s_3)\right]ds_3\right)^2 \right] \right]ds_1 \right)^2
dldx\]
Other terms can be derived analogously. They can be estimated in the similar way as in \cite{BGT3} though some new difficulties arise and the number of terms is substantially bigger. To obtain estimates we need the following inequalities
\begin{equation}
 \int_u^1 e^{-T(r-u)|z|^\alpha} \chi(r) dr \leq t-s, \quad 0\leq u \leq 1, \label{ineq:tightness1}
\end{equation}
\begin{equation}
 \intc{1}\int_u^1 e^{-T(r-u)|z|^\alpha} \chi(r) dr du \leq \frac{t-s}{T|z|^\alpha}, \label{ineq:tightness2}
\end{equation}
\begin{equation}
 \intc{u} e^{-T(u-s)|z|^\alpha} du \leq \frac{1}{T|z|^\alpha} \rbr{1-e^{-T|z|^\alpha}}, \label{ineq:tightness3}
\end{equation}
They are easily proved using inequality (\ref{tmp:chi_n}).\\
Now, to illustrate techniques required in estimations, we will carry out the proof for the term $R$ which is perhaps the most impressive one.\\
Firstly, we apply the Fubini theorem multiple times in order to separate the "time part" and the "space part"
\[R=\int_{0}^{+\infty} \int_{0}^{l} \int_{0}^{T} \int_{0}^{T} \int_{0}^{l} \int_{0}^{T} \int_{0}^{T} \chi_{T}(T - s_3)\chi_{T}(T - s_4)\chi_{T}(T - s_5)\chi_{T}(T - s_6)\:S\: ds_6ds_5ds_2ds_4ds_3ds_1dl \]
where
\[
S = \int_{\mathbb{R}^d}\mathcal{T}_{l-s_1}\left[\mathcal{T}_{s_1}\left[\mathcal{T}_{T-s_3}\left[\varphi_{T}(\cdot)\right]\right]\mathcal{T}_{s_1}\left[\mathcal{T}_{T-s_4}\left[\varphi_{T}(\cdot)\right]\right]\right]\mathcal{T}_{l-s_2}\left[\mathcal{T}_{s_2}\left[\mathcal{T}_{T-s_5}\left[\varphi_{T}(\cdot)\right]\right]\mathcal{T}_{s_2}\left[\mathcal{T}_{T-s_6}\left[\varphi_{T}(\cdot)\right]\right]\right]dx\]
By applying the Plancharel formula and definition (\ref{def: abbrev}) we get
\[
S = T^{-2} \int_{\mathbb{R}^{3d}}e^{-(l-s_1)|z|^\alpha}e^{-s_1|z_1|^\alpha}e^{-(T-s_3)|z_1|^\alpha}e^{-s_1|z -z_1|^\alpha}e^{-(T-s_4)|z -z_1|^\alpha}e^{-(l-s_2)|z|^\alpha}e^{-s_2|z_2|^\alpha}\]
\[
\times e^{-(T-s_5)|z_2|^\alpha}e^{-s_2|z -z_2|^\alpha}e^{-(T-s_6)|z -z_2|^\alpha}\widehat{\varphi}(z_1)\widehat{\varphi}(z -z_1)\widehat{\varphi}(z_2)\widehat{\varphi}(z -z_2)dz_2dz_1dz
\]
The Fubini theorem yields
\[R= T^{-2}\int_{\mathbb{R}^{3d}}\widehat{\varphi}(z_1)\widehat{\varphi}(z -z_1)\widehat{\varphi}(z_2)\widehat{\varphi}(z -z_2)\int_{0}^{+\infty} \int_{0}^{l} \int_{0}^{T} \int_{0}^{T} \int_{0}^{l} \int_{0}^{T} \int_{0}^{T} A ds_6ds_5ds_2ds_4ds_3ds_1dldz_2dz_1dz\]
where
\[
A = e^{-(l-s_1)|z|^\alpha}e^{-s_1|z_1|^\alpha}e^{-(T-s_3)|z_1|^\alpha}e^{-s_1|z -z_1|^\alpha}e^{-(T-s_4)|z -z_1|^\alpha}e^{-(l-s_2)|z|^\alpha}e^{-s_2|z_2|^\alpha}e^{-(T-s_5)|z_2|^\alpha}
\]
\[
\times e^{-s_2|z -z_2|^\alpha}e^{-(T-s_6)|z -z_2|^\alpha}\chi_{T}(T - s_3)\chi_{T}(T - s_4)\chi_{T}(T - s_5)\chi_{T}(T - s_6)
\]
Subsequent application of inequalities (\ref{ineq:tightness1}), (\ref{ineq:tightness3}), (\ref{ineq:tightness1}), (\ref{ineq:tightness3}) to integrals with respect to $s_6$, $s_5$, $s_4$, $s_3$ gives
\begin{equation}
R \leq (t-s)^2\int_{\mathbb{R}^{3d}}\widehat{\varphi}(z_1)\widehat{\varphi}(z -z_1)\widehat{\varphi}(z_2)\widehat{\varphi}(z -z_2)\frac{1}{|z_2|^\alpha}\frac{1}{|z_1|^\alpha} S(z,z_1, z_2)dz_2dz_1dz \label{ineq:R}
\end{equation}
where
\begin{equation}
S(z,z_1, z_2) = \int_{0}^{+\infty} \int_{0}^{l} \int_{0}^{l} e^{-(l-s_1)|z|^\alpha}e^{-s_1|z_1|^\alpha}e^{-s_1|z -z_1|^\alpha}e^{-(l-s_2)|z|^\alpha}e^{-s_2|z_2|^\alpha}e^{-s_2|z -z_2|^\alpha}ds_2ds_1dl, 
\end{equation}
A trivial verification shows that
\begin{eqnarray*}
S(z, z_1, z_2) = S_1(z,z_1, z_2) + S_2(z,z_1, z_2),\label{eq:S-decomposition}
\end{eqnarray*}
where
\[
S_1(z,z_1, z_2) = \frac{1}{2\norma{z}(\norma{z_1}+\norma{z-z_1}+\norma{z})(\norma{z_1}+\norma{z-z_1}+\norma{z_2}+\norma{z-z_2})}
\]
\[
S_2(z,z_1, z_2) = \frac{1}{2\norma{z}(\norma{z_2}+\norma{z-z_2}+\norma{z})(\norma{z_1}+\norma{z-z_1}+\norma{z_2}+\norma{z-z_2})}.
\]
Using (\ref{eq:S-decomposition}) we write the right-hand side of (\ref{ineq:R}) as $R_1 + R_2$, where $R_1$, $R_2$ have an obvious meaning. It is easy to see that
\[
R_1 = (t-s)^2\int_{\mathbb{R}^{3d}} \frac{\widehat{\varphi}(z_1)\widehat{\varphi}(z -z_1)}{\norma{z_1} |z-z_1|^{\alpha/2}} \frac{\widehat{\varphi}(z_2)\widehat{\varphi}(z -z_2)}{\norma{z_2}\norma{z-z_2}}\frac{1}{2|z|^{3/2\alpha}} dz_1 dz_2 dz
\]
Notice that function $f(x) = \frac{\widehat{\varphi}(x)}{\norma{x}}$ is square-integrable. The integral with respect to $z_2$ is is equal to $(f\ast f)(z)$. By Young's inequality (\ref{ineq:young}) it is easy to see that it is bounded (take $q_1 = q_2 = 2$), hence
\[
R_1 \leq c_1(t-s)^2\int_{\mathbb{R}^{d}} \frac{h(z)}{2|z|^{3/2\alpha}}  dz,
\]
where
\[
h(z) = \intr \frac{\widehat{\varphi}(z_1)\widehat{\varphi}(z -z_1)}{\norma{z_1} |z-z_1|^{\alpha/2}} dz_1 =  \left( \frac{\widehat{\varphi}(\cdot)}{\norma{\cdot}} \ast \frac{\widehat{\varphi}(\cdot)}{|\cdot|^{\alpha/2}} \right) (z)
\]
We may apply Young's inequality (\ref{ineq:young}) in two ways. Firstly taking $q_1 = 2/3$ and $q_2 = 3$ proves that $h$ is bounded, secondly taking $q_1 = q_2 = 1$ shows that $h$ is integrable. Hence 
\[
R_1 \leq c_2(t-s)^2.
\]
The proof for $R_2$ goes along the same lines.
\subsection{Proof of Theorem 2.2} \label{sec:calc3}
As the proof for the critical dimensions in the Poisson-starting system case is similar the one in Section \ref{sec:calc1}, we present only a sketch of the proof. Once again we follow the scheme described in Section \ref{sec:spacetime-schema}. (\ref{res:I1_critical}) can be obtained in the same way as \cite[(3.31)]{BGT3}. To prove the convergence (\ref{res:I2I3_critical}) of $I_2(T)$ one can follow the proof for the large dimension case and estimate arising terms $I_{211}$ and $I_{212}$ in a manner presented in \cite{BGT3} for $J_1$ and $J_2$ in the critical case. Limit for the $I_3$ is trivial.\\
Now we turn to the equilibrium case. We need to show (\ref{res:B_T_limit}).
\[
B_1(T) = \inti \intr \rbr{\T{t}\intc{T} \T{T-s_1} \varphi_T(x) \chi_T(T-s_1) ds_1 } \rbr{\T{t}\intc{T} \T{T-s_2} \varphi_T(x) \chi_T(T-s_2) ds_2 } dx dt
\]
\[
 = \inti \intc{T} \intc{T} \chi_T(T-s_1) \chi_T(T-s_2) \intr \T{t+T-s_1} \varphi_T(x) \T{t+T-s_2} \varphi_T(x) dx ds_1 ds_2 dt.
\]
By applying the Fourier transform we have
\[
B_1(T) = \frac{1}{(2\pi)^d} \inti \intc{T} \intc{T} \chi_T(T-s_1) \chi_T(T-s_2) \intr e^{(t+T-s_1)|z|^\alpha} e^{(t+T-s_2)|z|^\alpha} |\widehat{\varphi}_T(z)|^2 dz ds_1 ds_2 dt.
\]
Integrating with respect to $t$ yields
\begin{equation}
B_1(T) = c_1 \frac{A_T}{F^2_T} \label{res:limB_1},
\end{equation}
where
\[
A_T = \intr \frac{|\widehat{\varphi}(z)|^2}{|z|^\alpha} \rbr{\int_0^T e^{s|z|^\alpha} ds}^2 dz.
\]
Derivative of $A_T$ with respect to $T$ is given by
\[
A'_T = 2 \intr \frac{|\widehat{\varphi}(z)|^2}{|z|^\alpha} e^{-T|z|^\alpha} \frac{1}{|z|^\alpha} \rbr{1-e^{-T|z|^\alpha}} dz.
\]
In the critical case $\alpha = \frac{d}{2}$ so substituting $T^\frac{2}{d}z=z'$ we obtain
\[
A'_T = 2 \intr \frac{|\widehat{\varphi}(z'/T^\frac{2}{d})|^2}{|z'|^\alpha} e^{-|z'|^\alpha} \frac{1}{|z'|^\alpha} \rbr{1-e^{-|z'|^\alpha}} dz;
\]
 $\frac{1}{|z'|^\alpha} \rbr{1-e^{-|z'|^\alpha}}$ is bounded and $\frac{e^{-|z'|^\alpha}}{|z'|^\alpha}$ is integrable hence there exists a constant $c_2$
\[
A'_T \leq c_2.
\]
We obtain the limit of $B_1(T)$ using l'Hopital's rule ( $F^2_T = T \log T$)
\[
\lim_T B_1(T) = c_1 \lim_T \frac{A_T'}{\rbr{F^2_T}'} \leq \lim_T \frac{c_3}{\log T +1} = 0.
\]
This completes the proof of part 1). To show part 2) we, similarly as in the proof of Theorem \ref{thm: big-dimensions}, follow the scheme from Section \ref{sec:tightness}. In the same way we evaluate the terms arising from (\ref{ineq:tightness_main1}) and (\ref{ineq:tightness_main2}). Although the techniques to estimate them are similar to the ones presented in \cite{BGT3} we deal with more terms. We need the following estimates
\begin{equation}
\frac{1}{\log T} \intr \frac{f(z)}{|z|^{2\alpha}} \left[1- e^{-T|z|^\alpha} \right] dz \leq c(f) \label{ineq:tightness4},
\end{equation}
for $f$  bounded and integrable,
\begin{equation}
\frac{1}{\log T} \intr \frac{\widehat{\varphi}(z-z_1)}{|z-z_1|^\alpha} \left[1 - e^{-T|z-z_1|^\alpha} \right] \frac{1}{|z_1|^\alpha} \left[1 - e^{-T|z_1|^\alpha} \right] dz_1 \leq c(\varphi) \label{ineq:tightness5},
\end{equation}
for $\varphi$ rapidly decreasing,
\begin{equation}
\frac{1}{\log T} \intr \frac{\widehat{\varphi}(z-z_1)}{|z-z_1|^\alpha} \left[1 - e^{-T|z-z_1|^\alpha} \right] \frac{\widehat{\varphi}(z_1)}{|z_1|^\alpha} \left[1 - e^{-T|z_1|^\alpha} \right] dz_1 \leq f(z) \label{ineq:tightness6},
\end{equation}
where $f$ is integrable and bounded.\\
Inequalities (\ref{ineq:tightness4}) and (\ref{ineq:tightness5}) follows easily from l'Hopital's rule. To show (\ref{ineq:tightness6}) it suffices to observe that boundedness is a direct consequence of (\ref{ineq:tightness5}). The fact that $f\in \mathcal{L}^1$ follows from Young's inequality applied to 
\[
\intr \frac{\widehat{\varphi}(z-z_1)}{|z-z_1|^\alpha} \frac{\widehat{\varphi}(z_1)}{|z_1|^\alpha} dz_1.
\]
Finally, to illustrate problems arising in the critical dimension case, we show one example. Let us take the fourth term in $B^{IV}(0)$ (see \ref{res:BIV})
\[
\inti \intr V'''(0)(x,l) V'(0)(x,l) dx dl.
\]
One of terms resulting from its evaluation is
\begin{eqnarray*}
R = \int_{\mathbb{R}^d}\int_{0}^{+\infty} \mathcal{T}_{l}\left[\int_{0}^{T} \mathcal{T}_{T-s_1}\left[v'(0)(x, s_1)\int_{0}^{s_1} 
\mathcal{T}_{s_1-s_2}\left[v'(0)(x, s_2)v'(x, s_2)\right](x)ds_2\right](x)ds_1\right](x)\\
\mathcal{T}_{l}\left[v'(0)(x, T)\right](x)dldx.
\end{eqnarray*}
We substitute $v'(0)$ and change the order of integration
\[
R = \int_{0}^{+\infty} \int_{0}^{T} \int_{0}^{s_1} \int_{0}^{s_1} \int_{0}^{s_2} \int_{0}^{s_2} \int_{0}^{T} \chi_{T}(T - s_5)\chi_{T}(T - s_3)\chi_{T}(T - s_4)\chi_{T}(T - s_6) S ds_6ds_4ds_3ds_2ds_5ds_1dl,\]
where 
\[
S = \int_{\mathbb{R}^d}\mathcal{T}_{l}\left[\mathcal{T}_{T-s_1}\left[\mathcal{T}_{s_1-s_5}\left[\varphi_{T}(\cdot)\right]\mathcal{T}_{s_1-s_2}\left[\mathcal{T}_{s_2-s_3}\left[\varphi_{T}(\cdot)\right]\mathcal{T}_{s_2-s_4}\left[\varphi_{T}(\cdot)\right]\right]\right]\right]\mathcal{T}_{l}\left[\mathcal{T}_{T-s_6}\left[\varphi_{T}(\cdot)\right]\right]dx.
\]
By applying the Fourier transform we obtain 
\begin{eqnarray*}
S=\int_{\mathbb{R}^{3d}}e^{-l|z|^\alpha}e^{-(T-s_1)|z|^\alpha}e^{-(s_1-s_5)|z_1|^\alpha}e^{-(s_1-s_2)|z -z_1|^\alpha}e^{-(s_2-s_3)|z_2|^\alpha}e^{-(s_2-s_4)|z -z_1 -z_2|^\alpha}\\
\times e^{-l|z|^\alpha}e^{-(T-s_6)|z|^\alpha}\widehat{\varphi}(z_1)\widehat{\varphi}(z_2)\widehat{\varphi}(z -z_1 -z_2)\widehat{\varphi}(z)dz_2dz_1dz.
\end{eqnarray*} 
Once again we change the order of integration
\begin{equation}
R = T^{-2}\log T^{-2}\int_{\mathbb{R}^{3d}}\widehat{\varphi}(z_1)\widehat{\varphi}(z_2)\widehat{\varphi}(z -z_1 -z_2)\widehat{\varphi}(z) Q dz_2dz_1dz,  \label{res:uff}
\end{equation}
where
\begin{eqnarray*}
Q = \int_{0}^{+\infty} \int_{0}^{T} \int_{0}^{s_1} \int_{0}^{s_1} \int_{0}^{s_2} \int_{0}^{s_2} \int_{0}^{T} e^{-2l|z|^\alpha}e^{-(T-s_1)|z|^\alpha}e^{-(s_1-s_5)|z_1|^\alpha}
e^{-(s_1-s_2)|z -z_1|^\alpha}
e^{-(s_2-s_3)|z_2|^\alpha}\\
\times e^{-(s_2-s_4)|z -z_1 -z_2|^\alpha}
e^{-(T-s_6)|z|^\alpha}\chi_{T}(T - s_5)\chi_{T}(T - s_3)\chi_{T}(T - s_4)\chi_{T}(T - s_6)ds_6ds_4ds_3ds_2ds_5ds_1dl.
\end{eqnarray*}
By applying inequality (\ref{ineq:tightness1}) to the integral with respect $s_6$ we get
\begin{eqnarray*}
Q \leq c_1 T(t-s)\int_{0}^{+\infty} \int_{0}^{T} \int_{0}^{s_1} \int_{0}^{s_1} \int_{0}^{s_2} \int_{0}^{s_2} e^{-2l|z|^\alpha}e^{-(T-s_1)|z|^\alpha}e^{-(s_1-s_5)|z_1|^\alpha}e^{-(s_1-s_2)|z -z_1|^\alpha}\\
\times e^{-(s_2-s_3)|z_2|^\alpha}
e^{-(s_2-s_4)|z -z_1 -z_2|^\alpha}\chi_{T}(T - s_5)\chi_{T}(T - s_3)\chi_{T}(T - s_4)ds_4ds_3ds_2ds_5ds_1dl.
\end{eqnarray*}
Next we utilize (\ref{ineq:tightness3}) to eliminate the integral with respect to $s_5$
\begin{eqnarray*}
Q \leq c_2 T(t-s)\frac{1}{|z_1|^\alpha}\left[ 1 - e^{-T|z_1|^\alpha}\right]\int_{0}^{+\infty} \int_{0}^{T} \int_{0}^{s_1} \int_{0}^{s_2} \int_{0}^{s_2} e^{-2l|z|^\alpha}e^{-(T-s_1)|z|^\alpha}e^{-(s_1-s_2)|z -z_1|^\alpha}\\
\times e^{-(s_2-s_3)|z_2|^\alpha}e^{-(s_2-s_4)|z -z_1 -z_2|^\alpha}\chi_{T}(T - s_3)\chi_{T}(T - s_4)ds_4ds_3ds_2ds_1dl.
\end{eqnarray*}
Once again we use (\ref{ineq:tightness1}) this time to the integral with respect to $s_4$
\begin{eqnarray*}
Q \leq c_3 (t-s)^2 T^2 \frac{1}{|z_1|^\alpha}\left[ 1 - e^{-T|z_1|^\alpha}\right]\int_{0}^{+\infty} \int_{0}^{T} \int_{0}^{s_1} \int_{0}^{s_2} e^{-2l|z|^\alpha}e^{-(T-s_1)|z|^\alpha}e^{-(s_1-s_2)|z -z_1|^\alpha}\\
\times e^{-(s_2-s_3)|z_2|^\alpha}\chi_{T}(T - s_3)ds_3ds_2ds_1dl.
\end{eqnarray*}
Finally we apply (\ref{ineq:tightness3}) to the integrals with respect to $s_3$, $s_2$, $s_1$ consequently and integrate with respect to $l$
\begin{eqnarray*}
Q \leq c_4 (t-s)^2 T^{2}\frac{1}{|z_1|^\alpha}\left[ 1 - e^{-T|z_1|^\alpha}\right]\frac{1}{|z_2|^\alpha}\frac{1}{|z -z_1|^\alpha}\left[ 1 - e^{-T|z -z_1|^\alpha}\right]\frac{1}{|z|^{2\alpha}}\left[ 1 - e^{-T|z|^\alpha}.\right]
\end{eqnarray*}
We return to (\ref{res:uff})
\begin{eqnarray*}
R \leq c_5(t-s)^2\log T^{-2}\int_{\mathbb{R}^{3d}}\widehat{\varphi}(z_1)\widehat{\varphi}(z_2)\widehat{\varphi}(z -z_1 -z_2)\widehat{\varphi}(z)\frac{1}{|z_1|^\alpha}\left[ 1 - e^{-T|z_1|^\alpha}\right] \\
\times \frac{1}{|z_2|^\alpha}\frac{1}{|z -z_1|^\alpha}\left[ 1 - e^{-T|z -z_1|^\alpha}\right]\frac{1}{|z|^{2\alpha}}\left[ 1 - e^{-T|z|^\alpha}\right]dz_2dz_1dz.
\end{eqnarray*}
The integral with respect to $z_2$ is bounded
\[R\leq c_6 (t-s)^2\log T^{-2}\int_{\mathbb{R}^{2d}}\frac{\widehat{\varphi}(z_1)}{|z_1|^\alpha}\left[ 1 - e^{-T|z_1|^\alpha}\right]\frac{1}{|z -z_1|^\alpha}\left[ 1 - e^{-T|z -z_1|^\alpha}\right]\frac{\widehat{\varphi}(z)}{|z|^{2\alpha}}\left[ 1 - e^{-T|z|^\alpha}\right]dz_1dz.\]
Using inequality (\ref{ineq:tightness5}) we obtain 
\[R\leq c_7 (t-s)^2\log T^{-1}\int_{\mathbb{R}^{d}}\frac{\widehat{\varphi}(z)}{|z|^{2\alpha}}\left[ 1 - e^{-T|z|^\alpha}\right]dz.\]
We complete the proof by applying (\ref{ineq:tightness4}) and arriving at
\[R \leq c_8 (t-s)^2.\]


\begin{thebibliography}{}
\bibitem[1]{B}P. Billingsley, Convergence of Probability Measures., John Wiley\&Sons, New York, 1968.
\bibitem[2]{BZ1}M. Birkner and I. Z\"ahle, 
Functional central limit theorems for the occupation time of the origin for branching random walks in
$d\geq3$, Weierstra\ss{} Insitut f\"ur Angewandte Analysis und Stochastik,
Berlin, preprint No. 1011 (2005).
\bibitem[3]{BZ2}M. Birkner and I. Z\"ahle, 
A functional CLT for the occupation time of state-dependent branching random walk, to appear in Ann. Probab.
\bibitem[4]{BGR}T. Bojdecki, L.G. Gorostiza and S. Ramaswamy, Convergence of $\mathcal{S}'$-valued processes and space time random fields, J. Funct. Anal. 66 (1986), pp. 21-41.
\bibitem[5]{BGT1}T. Bojdecki, L.G. Gorostiza and A. Talarczyk, Sub-fractional Brownian motion and its relation to occupation times, Statist. Probab. Lett. 69 (2004), pp. 405-419.
\bibitem[6]{BGT2}T. Bojdecki, L.G. Gorostiza and A. Talarczyk, Limit theorems for occupation time fluctuations of branching systems I: Long-range dependence, Stoch. Proc. Appl. 116 (2006), pp. 1-18.
\bibitem[7]{BGT3}T. Bojdecki, L.G. Gorostiza and A. Talarczyk, Limit theorems for occupation time fluctuations of branching systems II: Critical and large dimensions Functional, Stoch. Proc. Appl. 116 (2006), pp. 19-35.
\bibitem[8]{BGT4}T. Bojdecki, L.G. Gorostiza and A. Talarczyk, A long range dependence stable process and an infinite variance branching system, 
www.arxiv.org, math.PR/0511739 (2005).
\bibitem[9]{BGT5}T. Bojdecki, L.G. Gorostiza and A. Talarczyk, Occupation time fluctuations of an infinite variance branching systems in large
dimensions, www.arxiv.org, math.PR/0511745 (2005).
\bibitem[10]{DW} J.D. Deuschel, K. Wang, Large deviations for the occupation time functional of a Poisson system of independent Brownian particles, Stoch. Proc. Appl. 52 (1994) 183-209.
\bibitem[11]{GW2}L.G. Gorostiza and A. Wakolbinger, 
Persistence criteria for a class of critical branching particle systems in continuous time.
Ann. Probab. 19 (1991), pp. 266-288.
\bibitem[12]{M}P. Miłoś, Occupation time fluctuations of Poisson and equilibrium finite variance branching systems, to appear in Probab. and Math. Stat.
\bibitem[13]{MI}I. Mitoma, Tightness of probabilities on $C\left(\left[0,1\right],\mathcal{S}'\right)$ and $D\left(\left[0,1\right],\mathcal{S}'\right)$, Ann. Probab. 11 (1983), pp. 989-999.
\end{thebibliography}
\end{document}